\documentclass[a4paper,12pt]{article}
\usepackage{amsmath}
\usepackage[english]{babel} 
\usepackage{graphicx,color} 
\DeclareGraphicsExtensions{.eps .pdf} 
\usepackage{epstopdf}
\usepackage{authblk}
\usepackage{float}
\usepackage{fixmath}
\usepackage{bm}
\usepackage{amsmath}
\usepackage{epstopdf}
\usepackage{float}
\usepackage{fixmath}
\usepackage{bm}
\usepackage{inputenc,t1enc}
\usepackage{bbm}
\usepackage{amsfonts}
\usepackage{latexsym}
\usepackage{mathrsfs}
\usepackage{inputenc,t1enc}
\usepackage{latexsym}
\DeclareSymbolFont{AMSb}{U}{msb}{m}{n}
\DeclareSymbolFontAlphabet{\mathbb}{AMSb}
\begin{document}  

\bibliographystyle{plain}

\newcommand{\nc}{\newcommand}
\nc{\nt}{\newtheorem}
\nt{defn}{Definition}
\nt{lem}{Lemma}
\nt{pr}{Proposition}
\nt{theorem}{Theorem}
\nt{cor}{Corollary}
\nt{ex}{Example}
\nt{ass}{Assumption}
\nt{step}{Step}
\nt{case}{Case}
\nt{subcase}{Subcase}
\nt{note}{Note}
\nc{\bd}{\begin{defn}} \nc{\ed}{\end{defn}}
\nc{\blem}{\begin{lem}} \nc{\elem}{\end{lem}}
\nc{\bpr}{\begin{pr}} \nc{\epr}{\end{pr}}
\nc{\bth}{\begin{theorem}} \nc{\eth}{\end{theorem}}
\nc{\bcor}{\begin{cor}} \nc{\ecor}{\end{cor}}
\nc{\bex}{\begin{ex}}  \nc{\eex}{\end{ex}}
\nc{\bass}{\begin{ass}}  \nc{\eass}{\end{ass}}
\nc{\bstep}{\begin{step}}  \nc{\estep}{\end{step}}
\nc{\bcase}{\begin{case}}  \nc{\ecase}{\end{case}}
\nc{\bsubcase}{\begin{subcase}}  \nc{\esubcase}{\end{subcase}}
\nc{\bnote}{\begin{note}}  \nc{\enote}{\end{note}}
\nc{\prf}{{\bf Proof.} }
\nc{\eop}{\hfill $\Box$ \\ \\}
\nc{\argmin}{\mathrm{argmin}}
\nc{\argmax}{\mathrm{argmax}}
\nc{\sgn}{\mathrm{sgn}}
\nc{\Var}{\mathrm{Var}}
\nc{\Cov}{\mathrm{Cov}}
\nc{\bak}{\!\!\!\!\!}
\nc{\IBD}{\mathrm{IBD}}
\nc{\supp}{\mathrm{supp}}
\nc{\dom}{\mathrm{dom}}
\nc{\R}{{\mathbb R}}
\nc{\peq}{\preceq}
\nc{\wt}{\widetilde}
\nc{\Mult}{\mathrm{Mult}}
\newcommand{\tcb}{\textcolor{black}}
\newcommand{\tcr}{\textcolor{black}}
\nc{\Prob}[1]{\mathbb{P}_{#1}}
\newcommand{\notex}[1]
{$^{(!)}$\marginpar[{\hfill\tiny{\sf{#1}}}]{\tiny{\sf{(!) #1}}}}
\newcommand\approxin
  {\raisebox{-1ex}{ $\stackrel{\textstyle\in}{\scriptstyle\sim}$ }}

\title{Estimating a monotone probability mass function with known flat regions}
\author[]{Dragi Anevski}
\author[]{Vladimir M. Pastukhov \thanks{pastuhov@maths.lth.se} }

\affil[]{Centre for Mathematical Sciences, Lund University,\\ Lund, Sweden}
\date{}
\maketitle
\begin{abstract}

We propose a new estimator of a discrete monotone probability mass function with known flat regions. We analyse its asymptotic properties and compare its performance to the Grenander estimator and to the monotone rearrangement estimator.
\end{abstract}

\section{Introduction}
In this paper we introduce a new estimator of a monotone discrete distribution. The problem has been studied before, and in particular by  \cite{jankowski2009estimation}, who were the first to study the estimation problem and who also introduced two new estimators. The problem of monotone probability mass function estimation is related to the problem of density estimation under shape constraints, first studied and much earlier by Grenander \cite{grenander}. The literature for the continuous case problem is vaste, to mention just a few results, see for example \cite{carolandyks, rao}.  In the discrete case problem some recent results are \cite{balabdur, durothuet, giguelay, jankowski2009estimation}.  Both in the discrete and continuous case problems one has derived, in particular, limit distribution results under the assumption of regions of constancy the true and underlying density/probability mass function. However, to our knowledge, one has previously not used the assumption of regions of constancy in the estimation procedure. In this paper, we do use this information in the constructing of the estimator. Thus we present a Maximum Likelihood estimator (MLE) under the assumption of regions of constancy of the probability mass function and derive some limit properties for new estimator.

The paper is mainly motivated by the paper by H. K. Jankowski and J. A. Wellner \cite{jankowski2009estimation}, which was the first to study the problem of estimating a discrete monotone distribution.

To introduce the estimator, suppose that $\bm{p} = \{ p_{i} \}_{i\in \mathbb{N}_{+}}$ is a monotone decreasing probability mass function with support $\mathbb{N}_{+}$ with several known flat regions, i.e. $p_{i} > 0$, $\sum_{i\in \mathbb{N}_{+}} p_{i} = 1$ and 
\begin{eqnarray}\label{restflat}
	p_{1} = \dots = p_{w_{1}} > p_{w_{1}+1} = \dots = p_{w_{1}+w_{2}} > \dots > p_{\sum_{j=1}^{m-1} w_{j} + 1} = \dots = p_{k},    
\end{eqnarray}
where $k = \sup\{ i: p_{i} > 0 \}$, $m$ is the number of flat regions of $\bm{p}$, $\bm{w} = (w_{1}, \dots, w_{m})$ is the vector of the lengths (the numbers of points) of the flat regions of the true mass function $\bm{p}$, so that $\sum_{j=1}^{m}w_{j}=k$ for $k < \infty$ and $\sum_{j=1}^{m}w_{j}=\infty$ otherwise. Note, that if $\bm{p}$ is strictly decreasing at some point $i$, then $w = 1$ and if $\bm{p}$ is strictly decreasing on the whole support, then $m=k$ and $\bm{w} = (1, \dots,1)$ for $k < \infty$ and $m=\infty$ otherwise.

Suppose that we have observed $X_{1}, X_{2}, \dots, X_{n}$ i.i.d. random variables with probability mass function $\bm{p}$. The empirical estimator of $\bm{p}$ is then given by
\begin{eqnarray}\label{unrMLEp}
	\hat{p}_{n, i} = \frac{n_{i}}{n}, \quad n_{i}= \sum_{j=1}^{n} 1 \{ X_{j} = i\}, \quad \text{$ i \in\mathbb{N}_{+}$,}
\end{eqnarray}
and it is also the unrestricted Maximum Likelihood Estimator (MLE)

\begin{eqnarray*}
\hat{\bm{p}}_{n} = \underset{\bm{g}\in\mathcal{G}_{k}}{\argmax} \, \prod_{i=1}^{k} g_{i}^{n_{i}},
\end{eqnarray*}
where 
\begin{eqnarray*}
\mathcal{G}_{k} = \Big\{  \bm{g}\in \mathbb{R}^{k}_{+}: \, \sum_{i=1}^{k}g_{i} = 1 \Big\}
\end{eqnarray*}
and with $n_{i} = \sum_{j=1}^{n} 1 \{ X_{j} = i\}$. Then for a given $n=\sum_{i=1}^{k}n_{i}$, the vector $(n_{1}, \dots, n_{k})$ follows a multinomial distributions Mult$(n, \bm{p})$.

The empirical estimator $\hat{\bm{p}}_{n}$ is unbiased, consistent and asymptotically 
normal, see \cite{jankowski2009estimation, van1998asymptotic}. It, however, does not guaranty that the order restriction
\begin{eqnarray}\label{ordres}
\hat{p}_{n,1} \geq \hat{p}_{n, 2} \geq \dots \geq \hat{p}_{n, k}
\end{eqnarray}
is satisfied.

We next discuss two estimators that do satisfy the order restrictions, first introduced in \cite{jankowski2009estimation}. These are  the order restricted MLE and monotone rearrangement of the empirical estimator.

The monotone rearrangement of the empirical estimator $\hat{\bm{p}}^{R}_{n}$ is defined as 
\begin{eqnarray}\label{rearp}
\hat{\bm{p}}^{R}_{n} = rear(\hat{\bm{p}}_{n}),
\end{eqnarray}
where $\hat{\bm{p}}_{n}$ is the unrestricted MLE in (\ref{unrMLEp}) and $rear(v)$ for a vector $v = (v_{1}, \dots, v_{k})$ is the reverse-ordered vector. The estimator  $\hat{\bm{p}}^{R}_{n}$ clearly satisfies the order restriction (\ref{ordres}).

The MLE  under the order restriction (\ref{ordres}), $\hat{\bm{p}}^{G}_{n}$,  is defined as
\begin{eqnarray*}
   \hat{\bm{p}}_{n}^G &=& \underset{\bm{f}\in\mathcal{F}_{k}}{\argmax} \, \prod_{i=1}^{k} f_{i}^{n_{i}},
\end{eqnarray*}
where 
\begin{eqnarray*}
	\mathcal{F}_{k} = \Big\{ \bm{f} \in \mathbb{R}^{k}_{+}: \, \sum_{i=1}^{k}f_{i} = 1, \, f_{1}  \geq f_{2} \geq \dots \geq f_{k} \Big\} \subset \mathcal{G}_{k} .
\end{eqnarray*}
It is equivalent to the isotonic regression  of the unrestricted MLE, see    \cite{barlowstatistical, jankowski2009estimation, robertsonorder}, defined by
\begin{eqnarray}\label{orMLEp}
\hat{\bm{p}}^{G}_{n} = \underset{\bm{f} \in \mathcal{F}_{k}}{\argmin}\sum_{i=1}^{k}[\hat{p}_{n, i} - f_{i}]^{2},
\end{eqnarray}
where  the basic estimator $\hat{p}_{n, i}$ is the unrestricted MLE in (\ref{unrMLEp}). The estimator $\hat{\bm{p}}^{G}_{n}$ is usually called the Grenander estimator and is derived using the same algorithm as for the continuous case problem: it is the vector of left derivatives of the least concave majorant (LCM) of the empirical distribution function $\mathbb{F}_{n}(x) = n^{-1}\sum_{i=1}^{n}1_{[1,x]}(X_{i})$.

The estimators $\hat{\bm{p}}^{G}_{n}$ and $\hat{\bm{p}}^{R}_{n}$ were introduced and studied in detail in the paper by Jankowski and Wellner \cite{jankowski2009estimation}. In particular, H. Jankowski and J. Wellner \cite{jankowski2009estimation} derived consistency of the estimators and analysed further asymptotic properties and performance of the estimators for different distributions and data sets. They showed  that $\sqrt{n}(\hat{\bm{p}}^{R}_{n} - \bm{p})$ and $\sqrt{n}(\hat{\bm{p}}^{G}_{n} - \bm{p})$ converge weakly to the processes $Y^{R}$ and $Y^{G}$ which are obtained by the following transform of a Gaussian process on the space $l_{2}$ with mean zero and covariance matrix with the components $p_{i}\delta_{ij} - p_{i}p_{j}$: for all periods of constancy, $r$ through $s$, of $\bm{p}$ let 
\begin{eqnarray*}
(Y^{R})^{(r,s)} &=& rear(Y^{(r,s)})\\
(Y^{G})^{(r,s)} &=& (Y^{(r,s)})^{G},
\end{eqnarray*}
where $Y^{(r,s)}$ denotes the $r$ through $s$ elements of $Y$, cf. Theorem 3.8 in \cite{jankowski2009estimation}.

In this paper we construct an estimator of a monotone probability mass function in the following way
\begin{eqnarray}\label{mainpr}
\hat{\bm{p}}^{*}_{n} = \underset{\bm{f}\in\mathcal{F}^{*}_{k}}{\argmax} \, \prod_{i=1}^{k} f_{i}^{n_{i}},
\end{eqnarray}
where 
\begin{eqnarray*}
	\mathcal{F}^{*}_{k} = \Big\{\bm{f} \in \mathbb{R}^{k}_{+}: \, &&\sum_{i=1}^{k}f_{i} = 1, \, f_{1} = \dots =  f_{w_{1}} \geq \\
	&&f_{w_{1}+1} =\dots = f_{w_{1}+w_{2}} \geq \dots \geq f_{\sum_{j=1}^{m-1} w_{j} + 1} = \dots = f_{k}\Big\}
\end{eqnarray*}
and with $n_{i} = \sum_{j=1}^{n} 1 \{ X_{j} = i\}$. We note that the vector $\bm{w} = (w_{1}, \dots, w_{m})$ constitutes the lengths of $m$ flat regions of the true probability mass function.

We propose the following algorithm:
\begin{enumerate}
\item Assume we are given a data set $(x_{1}, \dots, x_{n})$ of observations from $n$ i.i.d. random variables $X_{1}, X_{2}, \dots, X_{n}$ and the vector of the lengths of the flat regions of the true mass function $(w_{1}, \dots, w_{m})$.

\item We group the probabilities, which are required to be equal, at each flat region of $\bm{p} = \{ p_{i} \}_{i \in \{1, \dots, k \} }$ into the single parameters $p'_{j}$, $j \in \{1, \dots, m\}$. Note here, that the true values $p'_{j}$ are strictly decreasing and satisfy the following linear constraint
\begin{eqnarray*}
	w_{1}p'_{1}+w_{2}p'_{2} + \dots + w_{m}p'_{m} = 1.
\end{eqnarray*}

\item Next, we find the order restricted MLE of $\hat{\bm{p}}'_{n} = (p'_{1}, \dots p'_{m})$, which is equivalent to the isotonic regression with weights $(w_{1}, \dots, w_{m})$, 
\begin{eqnarray*}
\hat{\bm{p}}^{'G}_{n} = \underset{\bm{f}'\in\mathcal{F}'_{m,\bm{w}}}{\argmin}\sum_{j=1}^{m}[\hat{p}'_{n, j} - f'_{j}]^{2}w_{j},
\end{eqnarray*}
where 
\begin{eqnarray*}
	\mathcal{F}'_{m, \bm{w}} = \Big\{ \bm{f}' \in \mathbb{R}^{m}_{+}: \, \sum_{j=1}^{m}f'_{j}w_{j} = 1, \, f'_{1}  \geq f'_{2} \geq \dots \geq f'_{m} \Big\},
\end{eqnarray*}

and with $\hat{p}'_{n, j}$ the unrestricted MLE defined by
\begin{eqnarray*}\label{}
\hat{\bm{p}}'_{n} = \underset{\bm{g}' \in \mathcal{G}'_{m, \bm{w}}}{\argmax} \, \prod_{j=1}^{m}{g'}_{j}^{n'_{j}},
\end{eqnarray*}
where 
\begin{eqnarray*}\label{}
\mathcal{G}'_{m, \bm{w}}= \Big\{ \bm{g}' \in \mathbb{R}^{m}_{+}: \, \sum_{j=1}^{m}w_{j}g'_{j} = 1 \Big\} \supset \mathcal{F}'_{m, \bm{w}},
\end{eqnarray*}
cf. Lemma \ref{equivGren} below for a proof of the equivalence. Here the data are reduced to the vector $\bm{n}'$ with $n'_{j} = \sum_{l = 1}^{n} 1 \{ X_{l} \in \{ q_{j}, q_{j} + w_{j}\} \}$ and where $q_{j}$ is an index of the first element in the $j$-th flat region of $\bm{p}$.

\item Having obtained the MLE $\hat{\bm{p}}^{'G}_{n}$ of $\bm{p}'$, we finally construct the MLE $\hat{\bm{p}}^{*}_{n}$ of $\bm{p} \in \mathcal{F}^{*}_{k}$, by letting the probabilities in the flat region of $\bm{p}$ be equal to the corresponding values in $\hat{\bm{p}}^{'G}_{n}$. This can be written in matrix form as 
\begin{eqnarray}\label{phatstar}
\hat{\bm{p}}^{*}_{n} = \bm{A}\hat{\bm{p}}^{'G}_{n},
\end{eqnarray}
where $A$ is a $k \times m$ matrix, with non-zero elements all ones:
$$
[\bm{A}]_{q_{j}: q_{j}+w_{j} -1, j} = 1,
$$
with $j \in \{1, \dots, m \}$, $q_{j}$, $q_{j}$ is the first index of the $j$-th flat region of $\bm{p}$ and $w_{j}$ is the length of the $j$-th flat region.
\end{enumerate}

Our goal is to investigate the estimator $\hat{\bm{p}}^{*}_{n}$ and compare its performance with the monotone rearrangement estimator $\hat{\bm{p}}^{R}_{n}$ defined in (\ref{rearp}) and the Grenander estimator $\hat{\bm{p}}^{G}_{n}$ defined in (\ref{orMLEp}). 

The paper is organised as follows. In Lemma \ref{equivGren} in Section \ref{sec:2} we prove that the order restricted MLE, for the grouped parameters, is given by the isotonic regression of the unrestricted MLE of the grouped parameters. Next, Lemma \ref{asunMLE} shows consistency and asymptotic normality of the unrestricted MLE for the grouped parameters. After that in Lemma \ref{asorMLE} we show that the order restricted MLE for the grouped parameters is consistent and asymptotically Normal. Finally, in Theorem \ref{thm1} we show consistency and derive the limit distribution for the new estimator $\hat{\bm{p}}^{*}_{n}$. In Section \ref{sec:3} we make a comparison with previous estimators. In particular, in Lemma \ref{compdistl} we show that $\hat{\bm{p}}^{*}_{n}$ has, properly scaled, asymptotically smaller risk both with $l_2$ as well as with Hellinger loss, compared to the Grenander estimator. The asymptotically smaller risk of $\hat{\bm{p}}^{*}_{n}$ compared to $\hat{\bm{p}}^{R}_{n}$ follows from this result together with the result by   \cite{jankowski2009estimation} on the better risk performance of $\hat{\bm{p}}^{G}_{n}$  with respect to  $\hat{\bm{p}}^{R}_{n}$. The paper ends with a small simulation study, illustrating the small sample behaviour of $\hat{\bm{p}}^{*}_{n}$ in comparison with $\hat{\bm{p}}^{G}_{n}$ and $\hat{\bm{p}}^{R}_{n}$; the new estimator seems to perform better then both $\hat{\bm{p}}^{G}_{n}$ and $\hat{\bm{p}}^{R}_{n}$.

\section{Proof of characterization of estimator and asymptotic results}\label{sec:2}
In this section we prove the statements which have been made for the algorithm above and analyse the asymptotic properties of the estimator $\hat{\bm{p}}^{*}_{n}$. We begin with a Lemma which will be used later in this section.

\blem\label{propconv}
Assume $\{X_{n}\}$ and $\{Y_{n}\}$ are sequences of random variables, taking values in the metric space $(\mathbb{R}^{k}, l_{2})$ with $k\leq\infty$ endowed with its Borel sigma algebra. If $X_{n} \stackrel{d}{\to} X$ and $\lim_{n\to\infty}\mathbb{P}(X_{n} = Y_{n})=1$, then $Y_{n} \stackrel{d}{\to} X$.
\elem
\prf To prove the statement of the Lemma, we use the Portmanteau Lemma in \cite{van1998asymptotic}, giving several equivalent characterisations of distributional convergence. From the  Portmanteau Lemma it follows that we have to prove
$$
\mathbb{E}[h(Y_{n})] \to \mathbb{E}[h(X)]
$$
for all bounded Lipschitz functions $h$. By the triangle inequality
\begin{eqnarray}\label{pr11}
|\mathbb{E}[h(Y_{n})] - \mathbb{E}[h(X)]| \leq |\mathbb{E}[h(X_{n})] - \mathbb{E}[h(X)]| + |\mathbb{E}[h(Y_{n})] - \mathbb{E}[h(X_{n})]|,
\end{eqnarray}
where the first term $|\mathbb{E}[h(X_{n})] - \mathbb{E}[h(X)]| \to 0,$ by the Portmanteau Lemma.

Next, take an arbitrary $\varepsilon > 0$, then the second term in (\ref{pr11}) is bounded as
\begin{eqnarray}\label{pr12}
|\mathbb{E}[h(Y_{n})] - \mathbb{E}[h(X_{n})]| &\leq &\mathbb{E}[|h(Y_{n}) - h(X_{n})|]\nonumber \\
&\leq& \mathbb{E}[|h(Y_{n}) - h(X_{n})|1 \{ l_{2}(Y_{n}, X_{n}) > \varepsilon \}]\\ \nonumber
&&+  \mathbb{E}[|h(Y_{n}) - h(X_{n})|1 \{ l_{2}(Y_{n}, X_{n}) \leq \varepsilon \}].
\end{eqnarray}
Here using the boundness of $h$, for the first term in the right hand side of (\ref{pr12}) we have that
\begin{eqnarray*}\label{}
\mathbb{E}[|h(Y_{n}) - h(X_{n})|1 \{l_{2}(Y_{n}, X_{n}) > \varepsilon \}] &\leq& 2\sup \{h(X)\} \mathbb{E}[1 \{l_{2}(Y_{n}, X_{n}) > \varepsilon \}]\\
&=&  2\sup\{h(X)\} \mathbb{P}[l_{2}(Y_{n}, X_{n}) > \varepsilon ],
\end{eqnarray*}
where $\lim_{n\to\infty} \mathbb{P}[l_{2}(Y_{n}, X_{n}) > \varepsilon ] = 0$ for every $\varepsilon > 0$, since $\lim_{n\to\infty}\mathbb{P}(X_{n} = Y_{n})=1$.\newline

The second term in the right hand side of  (\ref{pr12}) can be written as 
\begin{eqnarray*}\label{}
\mathbb{E}[|h(Y_{n}) - h(X_{n})|1 \{l_{2}(Y_{n}, X_{n}) \leq \varepsilon \}] &\leq & \varepsilon |h|_{Lip} \mathbb{P}[l_{2}(Y_{n}, X_{n}) \leq\varepsilon ],
\end{eqnarray*}
where $|h|_{Lip}$ is the Lipschitz norm, i.e. $|h|_{Lip}$ is the smallest number such that $|h(x) - h(y)| \leq |h|_{Lip}l_{2}(x, y)$. Furthermore $\lim_{n\to\infty} \mathbb{P}[l_{2}(Y_{n}, X_{n}) \leq \varepsilon ] = 1$ for every $\varepsilon > 0$, since $\lim_{n\to\infty}\mathbb{P}(X_{n} = Y_{n})=1$.

Therefore, taking the limsup of the the left hand side of equation (\ref{pr12})  we obtain
$$
\underset{n\to \infty}{\limsup} |\mathbb{E}[h(Y_{n})] - \mathbb{E}[h(X_{n})]| \leq \varepsilon |h|_{Lip},
$$
where $\varepsilon$ is an arbitrary positive number. Thus
$$
|\mathbb{E}[h(Y_{n})] - \mathbb{E}[h(X_{n})]| \to 0,
$$
as $n\to\infty$.\eop

Our goal is to obtain the asymptotic distribution of $\hat{\bm{p}}^{*}_{n}$, defined in (\ref{mainpr}). The true probability mass function $\bm{p}$ satisfies the order restrictions in $\mathcal{F}^{*}_{k}$. Let us make a reparametrisation by grouping the probabilities, which are required to be equal, at each flat region of $\bm{p} = \{ p_{i} \}_{i \in \{1, \dots, k \} }$ into the single parameters $p'_{j}$, $j \in \{1, \dots, m\}$. The reparametrisation transforms $\mathcal{F}^{*}_{k}$ into
\begin{eqnarray*}
	\mathcal{F}'_{m, \bm{w}} = \Big\{ \bm{f}' \in \mathbb{R}^{m}_{+}: \, \sum_{j=1}^{m}f'_{j}w_{j} = 1, \, f'_{1}  \geq f'_{2} \geq \dots \geq f'_{m} \Big\}
\end{eqnarray*}
and the the estimation problem (\ref{mainpr}) becomes
\begin{eqnarray}
	\hat{\bm{p}}^{'G}_{n} = \underset{\bm{f}'\in\mathcal{F}'_{m, \bm{w}}}{\argmax} \, { f'}_{1}^{n'_{1}}{f'}_{2}^{n'_{2}}\cdots {f'}_{m}^{n'_{m}},\label{repMLE}
\end{eqnarray}
where $n'_{j} = \sum_{l = 1}^{n} 1 \{ X_{l} \in \{ q_{j}, q_{j} + w_{j} - 1\} \}$ with $q_{j}$ an index of the first element in the $j$-th flat region of $\bm{p}$.

\blem \label{equivGren}
The solution $\hat{\bm{p}}^{'G}_{n}$ to the ML problem, defined in (\ref{repMLE}), is given by the weighted isotonic regression problem
\begin{eqnarray}\label{orMLEppr}
	\hat{\bm{p}}^{'G}_{n} = \underset{\bm{f}' \in \mathcal{F}'_{m, \bm{w}}}{\argmin}\sum_{j=1}^{m}[\hat{p}'_{n, j} - f'_{j}]^{2}w_{j},
\end{eqnarray}
where $\hat{p}'_{n, j}$ is the unrestricted (without order restrictions) MLE
\begin{eqnarray}\label{unrMLEppr}
	\hat{\bm{p}}'_{n} = \underset{\bm{g}' \in \mathcal{G}'_{m, \bm{w}}}{\argmax} \, { g'}_{1}^{n'_{1}}{g'}_{2}^{n'_{2}}\cdots {g'}_{m}^{n'_{m}},
\end{eqnarray}
where 
\begin{eqnarray*}\label{}
\mathcal{G}'_{m, \bm{w}} = \Big\{ \bm{g}' \in \mathbb{R}^{m}_{+}: \, \sum_{j=1}^{m}w_{j}g'_{j} = 1 \Big\}.
\end{eqnarray*}
\elem
\prf The result is the consequence of the problem of maximising the product of several factors, given relations of order and linear side condition, cf. pages 45--46 in \cite{barlowstatistical} and pages 38--39 \cite{robertsonorder}. In fact the results show, that for a product of several factors the MLE under the order restrictions coincides with the isotonic regression of the unrestricted ML estimates.

\eop

Next, we analyse the asymptotic behaviour of the unrestricted MLE $\hat{\bm{p}}'_{n}$ in (\ref{unrMLEppr}).

\blem \label{asunMLE}
The unrestricted MLE $\hat{\bm{p}}'_{n}$ in (\ref{unrMLEppr}) is given by 
$$
\hat{p}'_{n, j} = \frac{n'_{j}}{w_{j}n},
$$
$n'_{j} = \sum_{l = 1}^{n} 1 \{ X_{l} \in \{ q_{j}, q_{j} + w_{j} - 1\} \}$, where $q_{j}$ is an index of the first element in the flat region of $\bm{p}$.\\
It is consistent 
$$
\hat{\bm{p}}'_{n} \stackrel{p}{\to} \bm{p'} 
$$
and asymptotically normal
$$
\sqrt{n}(\hat{\bm{p}}'_{n} - \bm{p'}) \stackrel{d}{\to} \mathcal{N}(\bm{0},  \Sigma' ),
$$
where $\Sigma' $ is an $m \times m$ matrix such that $[\Sigma']_{ij} = \delta_{ij}\frac{p'_{i}}{w_{i}} - p'_{i}p'_{j}$, with $\delta_{ij}$ the indicator function for $i=j$.
\elem

\prf 
The result of the Lemma  for a case of a finite support of $\bm{p}$ ($k < \infty$ and consequently $m < \infty$)  follows directly from the Theorem 2 in  \cite{aitchison1958maximum}, also see pages 79--82 in \cite{silvey1975statistical}. 

Next, we consider a case of an infinite support of  $\bm{p}$ ($k=\infty$ and, obviously, $m=\infty$).  Let us introduce the notations $Z_{n} = \sqrt{n}(\hat{\bm{p}}'_{n} - \bm{p'})$ and $Z $ for a    $\mathcal{N}(\bm{0},  \Sigma')$-distributed r.v. and note that $\{Z_{n} \}$ is a sequence of processes in $l_{2}$, endowed with its Borel sigma algebra $\mathscr{B}$. 

First, for any finite integer $s$ the sequence of vectors $(Z_{n,t})_{t=1}^{s}$ converges in distribution to the vector $(Z)_{t=1}^{s} \in  \mathcal{N}(\bm{0},  \Sigma'_{s})$, where $[\Sigma'_{s}]_{ij} = \delta_{ij}\frac{p_{i}}{w_{i}} - p_{i}p_{j}$ with $i = 1, \dots, s$ and $j = 1, \dots, s$. This fact follows again from  \cite{aitchison1958maximum, silvey1975statistical}.

Second, we show that the sequence $\{ Z_{n} \}$ is tight in the $l_{2}$-norm metric. This is shown similarly to as in \cite{jankowski2009estimation}. In fact, from Lemma 6.2 in \cite{jankowski2009estimation} it is  enough to show that the two conditions 
\begin{eqnarray*}
\sup_{n}\mathbb{E}[|| Z_{n}||_{2}^{2}] &<& \infty, \\
\lim_{r\to\infty}\sup_{n}\sum_{t\geq r}\mathbb{E}[|Z_{n,t}|^2] &= &0,
\end{eqnarray*}
are satisfied. We note, that for any $n$ 
\begin{eqnarray*}
Z_{n,j}& =& \sqrt{n}(\hat{p}'_{n, j} - p_{j}') \\
&= &\sqrt{n}(\frac{n'_{j}}{w_{j}n} - p_{j}'),
\end{eqnarray*}
where $n'_{j}$  is Bin$(n, w_{j}p'_{j})$-distributed. Therefore, $\mathbb{E}[Z_{n,j}^2] = \frac{p'_{j}}{w_{j}} - p'_{j}p'_{j}$. Thus, both conditions of Lemma 6.2 are satisfied.

Third, since the space $l_{2}$ is separable and complete, from Prokhorov's theorem \cite{shiryaev} it follows that $\{Z_{n} \}$ is relatively compact, which means that every sequence from $\{Z_{n} \}$ contains a subsequence, which converges weakly to some process $Z$. In addition, if the limit processes have the same laws for every convergent subsequence, then $\{Z_{n} \}$  converges weakly to $\{Z \}$. 

 Next, we show the equality of laws of the limit processes of the convergent subsequences. First, note that since $l_{2}$ is a separable space,  the Borel $\sigma$-algebra equals the $\sigma$-algebra generated by open balls in $l_{2}$ \cite{bogachev}. Then, it is enough to show that the limit laws agree on finite intersections of open balls in $l_{2}$, since these constitute a  $\pi$-system. To show this, we note that the open balls in $l_{2}$ can be written as
\begin{eqnarray*}
B(z, \varepsilon) &=&  \cap_{M\geq 1}B_M,
\end{eqnarray*}
where 
\begin{eqnarray*}
   B_M&=&\cup_{n \geq 1}A_n^{M}, \\
   A_n^{M}&=&\{y \in l_{2}: \sum_{j=1}^{M}|z_{j} - y_{j}|^{2} < \varepsilon - \frac{1}{n}  \}.
\end{eqnarray*}

By the finite support part of the Lemma, the vectors $Z_n^{(M)}=(Z_{n,t})_{t=1}^{M}$ converge weakly to $Z^{(M)}=(Z_t)_{t=1}^M$ for all finite $M$, which implies that any subsequence of $Z_n^{(M)}$ converges weakly to $Z^{(M)}$. That means that, with $P_{n}^{(M)}$ the law of an arbitrary but fixed subsequence of $Z_n$, and $P^{(M)}$ the law of $Z^{(M)}$, $P_{n}^{(M)}(A)\to P^{(M)}(A)$ for any $P^{(M)}$-continuity set $A$. We note that the limit law $P^{(M)}$ is the same for all subsequences. Therefore, since $A_n^{M}$ is a continuity set for the Gaussian limit law $P^{(M)}$, and by the continuity properties of a probability measure, we obtain
\begin{eqnarray*}
      P(B(z, \varepsilon))&=&\lim_{M\to \infty}P(B_M)\\
            &=&\lim_{M\to \infty} \lim_{n \to \infty} P(A_n^M)\\
      &=&\lim_{M\to \infty} \lim_{n \to \infty} P^{(M)}(A_n^M),
\end{eqnarray*}
where $P$ is the law of $Z$. 
Thus, we have shown that the limit laws, $P$, of the convergent subsequences of $\{ Z_{n} \}$ agree on the open balls $B(z, \varepsilon)$, and, therefore, also on the finite intersections of these open balls. Since the laws agree on the $\pi$-system (they are all equal to $P$), they agree on the Borel $\sigma$-algebra. \eop

Summarising the results from the previous Lemmas, we obtain the final limit result for the estimator $\hat{\bm{p}}^{'G}_{n}$.
\blem \label{asorMLE}
The estimator $\hat{\bm{p}}^{'G}_{n}$ is consistent 
$$
\hat{\bm{p}}^{'G}_{n} \stackrel{p}{\to} \bm{p'} 
$$
and asymptotically normal
$$
\sqrt{n}(\hat{\bm{p}}^{'G}_{n} - \bm{p'}) \stackrel{d}{\to} \mathcal{N}\Bigg(\bm{0},  \Sigma \Bigg),
$$
where $[\Sigma]_{ij} = \delta_{ij}\frac{p_{i}}{w_{i}} - p_{i}p_{j}$ and $\delta_{ij}$ is the indicator function for $i=j$.
\elem
\prf From Lemma \ref{asunMLE} it follows that the basic estimator $\hat{\bm{p}}'_{n}$ is consistent
$$
\hat{\bm{p}}'_{n} \stackrel{p}{\to} \bm{p'}.
$$

From Theorem 2.2 in \cite{barlowstatistical}, it follows that if the basic estimator is consistent, then its isotonic regression is also consistent 
$$
\hat{\bm{p}}^{'G}_{n} \stackrel{p}{\to} \bm{p'}. 
$$

Since $\hat{\bm{p}}'_{n}$ and $\hat{\bm{p}}^{'G}_{n}$ both are consistent and since $\bm{p'}$ is an interior point of $ \mathcal{F}'_{m, \bm{w}} \subset \mathcal{G}'_{m, \bm{w}}$, there is an open set $\omega \subset \mathcal{F}'_{m, \bm{w}}$ such that $\bm{p'} \in \omega$ and 
\begin{eqnarray*}\label{}
\mathbb{P}[\hat{\bm{p}}'_{n} \in \omega] &\to& 1,\\
\mathbb{P}[\hat{\bm{p}}^{'G}_{n} \in \omega] &\to& 1,
\end{eqnarray*}
as $n\to\infty$. Furthermore, since $\omega \subset \mathcal{F}'_{m, \bm{w}}$ and  as long as $\hat{\bm{p}}'_{n} \in \mathcal{F}'_{m, \bm{w}}$, the equality $\hat{\bm{p}}^{'G}_{n} = \hat{\bm{p}}'_{n}$ holds, we have that
\begin{eqnarray*}\label{}
\mathbb{P}[\hat{\bm{p}}'_{n} \in \omega, \hat{\bm{p}}^{'G}_{n} \in \omega] = \mathbb{P}[\hat{\bm{p}}'_{n} \in \omega, \hat{\bm{p}}^{'G}_{n} \in \omega,\hat{\bm{p}}'_{n} = \hat{\bm{p}}^{'G}_{n} ] \leq \mathbb{P}[\hat{\bm{p}}'_{n} = \hat{\bm{p}}^{'G}_{n} ]
\end{eqnarray*}
and, since the left hand side of this inequality goes to one as $n\to\infty$, we have shown that $\mathbb{P}_{n \to \infty}[\hat{\bm{p}}^{'G}_{n} = \hat{\bm{p}}'_{n}] = 1$.

Now let $X_{n} = \sqrt{n}(\hat{\bm{p}}'_{n} - \bm{p'})$ and $Y_{n} = \sqrt{n}(\hat{\bm{p}}^{'G}_{n} - \bm{p'})$. Then, clearly,
\begin{eqnarray*}\label{}
\mathbb{P}[X_{n} = Y_{n}] = \mathbb{P}[\hat{\bm{p}}^{'G}_{n}  = \hat{\bm{p}}'_{n}]\to 1,
\end{eqnarray*}
as $n\to\infty$. Applying Lemma \ref{propconv} shows the statement of the Lemma. \eop

\bth\label{thm1}
The estimator $\hat{\bm{p}}^{*}_{n} = \bm{A}\hat{\bm{p}}^{'G}_{n}$ is consistent 
$$
\hat{\bm{p}}^{*}_{n} \stackrel{p}{\to} \bm{p} 
$$
and asymptotically normal
$$
\sqrt{n}(\hat{\bm{p}}^{*}_{n} - \bm{p}) \stackrel{d}{\to} \mathcal{N}(\bm{0},  \Sigma^{*} ),
$$
where $\Sigma^{*} = A\Sigma A^{T}$ with $[\Sigma]_{ij} = \delta_{ij}\frac{p_{i}}{w_{i}} - p_{i}p_{j}$ and $\bm{A}$ is a $k\times m$ matrix whose non-zero elements are $[\bm{A}]_{q_{j}: q_{j}+w_{j} -1, j} = 1$, $j \in \{1, \dots, m \}$, $q_{j}$ is the first index of the $j$-th flat region of the true mass function $\bm{p}$ and $w_{j}$ stands for the $j$-th regions length.
\eth
\prf From Lemma \ref{asorMLE} it follows that  $\hat{\bm{p}}'_{n}$ is consistent and asymptotically normal. The estimator $\hat{\bm{p}}^{*}_{n}$ is given by (\ref{phatstar}).
The statements of the Theorem now follow from the Delta method (see, for example, Theorem 3.1 in \cite{van1998asymptotic}). \eop

\section{Comparison of the estimators}\label{sec:3}
To compare the estimators we consider the $l_{2}$ metric
$$
l_{2}^{2}(\hat{\bm{p}}, \bm{p}) = ||\hat{\bm{p}} - \bm{p}||_{2}^{2} =\sum_{i=1}^{k}(\hat{p}_{i} - p_{i})^2,
$$
with $k\leq \infty$, and the Hellinger distance
$$
H^{2}(\hat{\bm{p}}, \bm{p}) = \frac{1}{2}\sum_{i=1}^{k}(\sqrt{\hat{p}_{i}} - \sqrt{p_{i}})^2,
$$
with $k < \infty$.  
In \cite{jankowski2009estimation} it has been shown that the Grenander estimator $\bm{p}^{G}_{n}$ has smaller risk than the rearrangement estimator $\bm{p}^{R}_{n}$, for both $l_{2}$ and $H$ loss. 

The next lemma shows that the new estimator $\hat{\bm{p}}^{*}_{n}$ performs better than the Grenander estimator $\bm{p}^{G}_{n}$, asymptotically, in both the expected $l^2$ and Hellinger distance sense, properly normalised.
\blem \label{compdistl}
For the $l_{2}$ metric we have that\\
\begin{eqnarray*}
\lim_{n\to\infty} \mathbb{E}[n l_{2}^{2}(\hat{\bm{p}}^{*}_{n}, \bm{p})]  \leq \lim_{n\to\infty} \mathbb{E}[n l_{2}^{2}(\hat{\bm{p}}^{G}_{n}, \bm{p})]
\end{eqnarray*}
and for the Hellinger distance with $k < \infty$,
\begin{eqnarray*}
\lim_{n\to\infty} \mathbb{E}[n H^{2}(\hat{\bm{p}}^{*}_{n}, \bm{p})]  \leq \lim_{n\to\infty} \mathbb{E}[n H^{2}(\hat{\bm{p}}^{G}_{n}, \bm{p})] 
\end{eqnarray*}
Equalities hold if and only if the true probability mass function $\bm{p}$ is strictly monotone.
\elem
\prf First, from Theorem \ref{thm1} and the continuous mapping theorem we have
\begin{eqnarray}
n l_{2}^{2}(\hat{\bm{p}}^{*}_{n}, \bm{p}) \stackrel{d}{\to}  ||V||_{2}^{2},\label{l2pst}
\end{eqnarray}
where $V \in \mathcal{N}(\bm{0},  \Sigma^{*} )$.

Second, using the reduction of error property of isotonic regression (Theorem 1.6.1 in \cite{robertsonorder}), for any $n$ we have
$$
\sum_{j=1}^{m}[\hat{p}^{'G}_{n, j} - p'_{j}]^{2}w_{j} \leq \sum_{j=1}^{m}[\hat{p}'_{n, j} - p'_{j}]^{2}w_{j},
$$
which is the same as
$$
n l_{2}^{2}(\hat{\bm{p}}^{*}_{n}, \bm{p}) \leq n l_{2}^{2}(\bar{\bm{p}}_{n}, \bm{p}),
$$
where $\bar{\bm{p}}_{n} = \bm{A}\hat{\bm{p}}'_{n}$ is constructed from $\hat{\bm{p}}'_{n}$ in the same way as $\hat{\bm{p}}^{*}_{n}$ from $\hat{\bm{p}}^{'G}_{n}$. Since that for every $M>0$ we have
\begin{eqnarray*}
\mathbb{E}[n l_{2}^{2}(\hat{\bm{p}}^{*}_{n}, \bm{p})1 \{ n l_{2}^{2}(\hat{\bm{p}}^{*}_{n}, \bm{p}) > M\}] &\leq& \mathbb{E}[n l_{2}^{2}(\bar{\bm{p}}_{n}, \bm{p})1 \{ n l_{2}^{2}(\hat{\bm{p}}^{*}_{n}, \bm{p}) > M\}] \leq\\
&& \mathbb{E}[n l_{2}^{2}(\bar{\bm{p}}_{n}, \bm{p})1 \{ n l_{2}^{2}(\bar{\bm{p}}_{n}, \bm{p}) > M\}]
\end{eqnarray*}
and 
\begin{eqnarray}
&&\limsup_{n\to\infty} \mathbb{E}[n l_{2}^{2}(\hat{\bm{p}}^{*}_{n}, \bm{p})1 \{ n l_{2}^{2}(\hat{\bm{p}}^{*}_{n}, \bm{p}) > M\}] \leq \label{lmspl2} \\
&&\limsup_{n\to\infty} \mathbb{E}[n l_{2}^{2}(\bar{\bm{p}}_{n}, \bm{p})1 \{ n l_{2}^{2}(\bar{\bm{p}}_{n}, \bm{p}) > M\}].\nonumber
\end{eqnarray}
From Lemma \ref{asunMLE} for any $n$ we have that
$$
\mathbb{E}[n l_{2}^{2}(\bar{\bm{p}}_{n}, \bm{p})] =  \sum_{j=1}^{m} w_{j}p_{j}(\frac{1}{w_{j}} - p_{j}) < \infty,
$$
and also
$$
\mathbb{E}[||V||_{2}^{2}] = \sum_{j=1}^{m} w_{j}p_{j}(\frac{1}{w_{j}} - p_{j})
$$
and using the Delta method and the continuous mapping theorem \cite{van1998asymptotic} it can be shown that $n l_{2}^{2}(\bar{\bm{p}}_{n}, \bm{p}) \stackrel{d}{\to}  ||V||_{2}^{2}$
which proves that the sequence $\{n l_{2}^{2}(\bar{\bm{p}}_{n}, \bm{p})\}$ is  asymptotically uniformly integrable (see, for example, Theorem 2.20 in \cite{van1998asymptotic})
$$
\lim_{M\to\infty} \limsup_{n\to\infty} \mathbb{E}[n l_{2}^{2}(\bar{\bm{p}}_{n}, \bm{p})1 \{ n l_{2}^{2}(\bar{\bm{p}}_{n}, \bm{p}) > M\}] = 0,
$$ 
which together with (\ref{lmspl2}) proves that 
$$
\lim_{M\to\infty}\limsup_{n\to\infty} \mathbb{E}[n l_{2}^{2}(\hat{\bm{p}}^{*}_{n}, \bm{p})1 \{ n l_{2}^{2}(\hat{\bm{p}}^{*}_{n}, \bm{p}) > M\}] = 0,
$$
which shows the asymptotic uniform integrability of the sequence $n l_{2}^{2}(\hat{\bm{p}}^{*}_{n}, \bm{p})$.

Third, since the sequence ${n l_{2}^{2}(\hat{\bm{p}}^{*}_{n}, \bm{p})}$ is asymptotically uniformly integrable and converges in distribution to $V$, it also converges in expectation (Theorem 2.20 in \cite{van1998asymptotic})
$$
\lim_{n\to\infty} \mathbb{E}[n l_{2}^{2}(\hat{\bm{p}}^{*}_{n}, \bm{p})]  =  \mathbb{E}[||V||_{2}^{2}] = \sum_{j=1}^{m} w_{j}p'_{j}(\frac{1}{w_{j}} - p'_{j})
$$
Furthermore, for $\hat{\bm{p}}^{G}_{n}$, proposed in \cite{jankowski2009estimation}, we have 
$$
\lim_{n\to\infty} \mathbb{E}[n l_{2}^{2}(\hat{\bm{p}}^{G}_{n}, \bm{p})] = \sum_{j=1}^{m} \sum_{q=1}^{w_{j}} p'_{j}(\frac{1}{q} - p'_{j}).
$$
It is obvious that $\sum_{j=1}^{m} w_{j}p'_{j}(\frac{1}{w_{j}} - p'_{j}) \leq \sum_{j=1}^{m} \sum_{q=1}^{w_{j}} p'_{j}(\frac{1}{q} - p'_{j})$. This finishes the proof of statement for the $l_{2}$ metric.

To prove the  statement for Hellinger distance, let us assume that $k < \infty$ is an arbitrary. It is sufficient to note that  
$$
nH^{2}(\hat{\bm{p}}^{*}_{n}, \bm{p}) = \frac{1}{2} \sum_{i=1}^{k}\frac{[\sqrt{n}(p^{*}_{n,i} - p_{i})]^{2}}{(\sqrt{p^{*}_{n,i}} + \sqrt{p_{i}})^{2}},
$$
since then, from the weak convergence and consistency of $\hat{\bm{p}}_n^{*}$, Slutsky's theorem and the continuous mapping theorem, it follows that 
$$
nH^{2}(\hat{\bm{p}}^{*}_{n}, \bm{p}) \stackrel{d}{\to} \frac{1}{8}\sum_{i=1}^{k}\frac{V_{i}^{2}}{p_{i}}.
$$

Furthermore, asymptotic uniform integrability of $nH^{2}(\hat{\bm{p}}^{*}_{n}, \bm{p})$ can be shown using the inequality $nH^{2}(\hat{\bm{p}}^{*}_{n}, \bm{p}) \leq n l_{1}^{2}(\hat{\bm{p}}^{*}_{n}, \bm{p})$, and asymptotic integrability of $n l_{1}^{2}(\hat{\bm{p}}^{*}_{n}, \bm{p})$, see \cite{jankowski2009estimation}. Therefore we also have convergence in expectation
\begin{eqnarray}
\lim_{n\to\infty} \mathbb{E}[n H^{2}(\hat{\bm{p}}^{*}_{n}, \bm{p})] = \mathbb{E}\Big[\sum_{i=1}^{k}\frac{1}{8}\frac{V_{i}^{2}}{p_{i}} \Big] = \frac{1}{8}\sum_{j=1}^{m}w_{j}(\frac{1}{w_{j}} - p_{j}). \label{eq:H-norm}
\end{eqnarray}

Finally, \cite{jankowski2009estimation} shows that the Hellinger distance of the estimator $\hat{\bm{p}}^{G}_{n}$ converges in expectation
$$
\lim_{n\to\infty} \mathbb{E}[n H^{2}(\hat{\bm{p}}^{G}_{n}, \bm{p})] = \frac{1}{8}\sum_{j=1}^{m}\sum_{q=1}^{w_j}(\frac{1}{q} - p_{j}) \geq \lim_{n\to\infty} \mathbb{E}[n H^{2}(\hat{\bm{p}}^{*}_{n}, \bm{p})]
$$
where we note the inequality from a comparison with $(\ref{eq:H-norm})$. It is clear that equality holds if and only if $\bm{p}$ is strictly monotone. \eop

\begin{figure}[t]
\center
\includegraphics[scale=0.8]{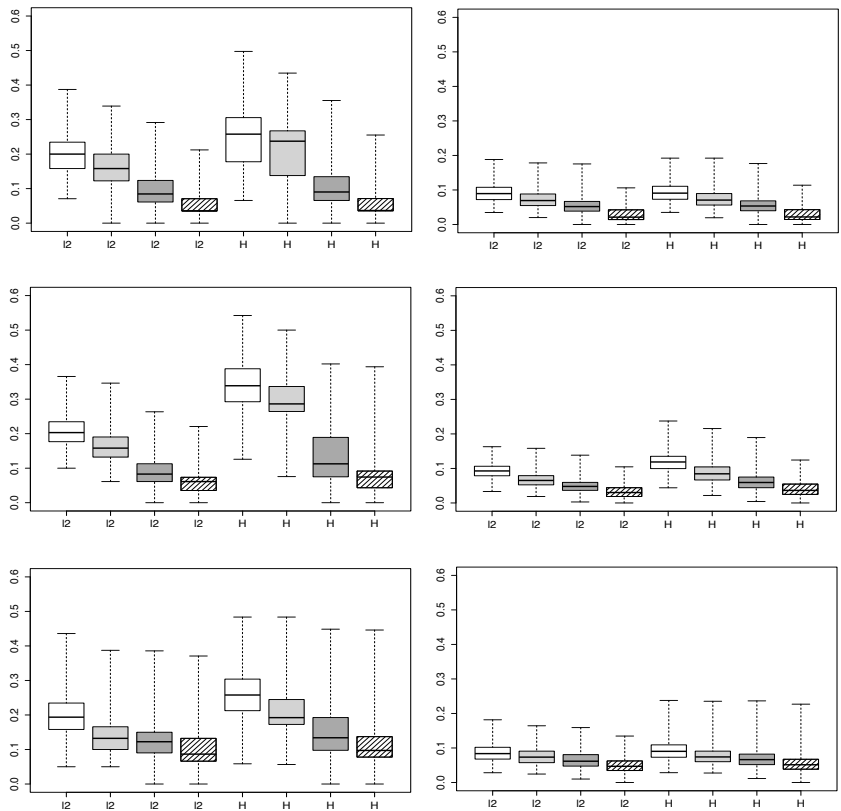}
\caption{The boxplots for $l_2$ norms and Hellinger distances for the estimators: the empirical estimator $\hat{\bm{p}}'$ (white), the rearrangement estimator $\hat{\bm{p}}^{R}_{n}$ (grey), Grenander estimator $\hat{\bm{p}}^{G}_{n}$ (dark grey) and estimator $\hat{\bm{p}}^{*}_{n}$ (shaded).}
\label{estims}
\end{figure}

For a visualisation of the finite sample performance of the proposed estimator $\hat{\bm{p}}_n^{*}$, we make a small simulation study. We choose the same probability mass functions as the ones chosen in \cite{jankowski2009estimation}. In Figure \ref{estims} we present results of Monte Carlo simulations for 1000 samples, for sample sizes $n=20$ and $n=100$, for the probability mass functions
\begin{enumerate}
\item (top) $p(x) = 0.2 U(4) + 0.8 U(8)$, 
\item (center) $p(x) = 0.15 U(4) + 0.1 U(8) + 0.75 U(12)$, 
\item (bottom) $p(x) = 0.25 U(2) + 0.2 U(4) + 0.15 U(6) + 0.4 U(8)$,
\end{enumerate}
where $U(k)$ stands for the uniform discrete distribution on $\{1, \dots, k\}$. The results shown are  boxplots for the Hellinger distance and $l_2$ metric, with sample sizes $n=20$ on the left and $n=100$ on the right in Fig. \ref{estims}.

The simulation study clearly illustrates that the newly proposed estimator $\hat{\bm{p}}^{*}_{n}$ has a better finite sample performance than both the Grenander and the monotone rearrangement estimators, in both $l_{2}$ and $H$ distance sense.

\section*{Acknowledgements}
VP's research is fully supported and  DA's research is partially supported by the Swedish Research Council, whose support is gratefully acknowledged.

\end{document}